\documentclass[12pt,psamsfonts]{amsart}
\usepackage{amsmath,amsthm,amsfonts,amssymb}
\usepackage{eucal}

\newcommand{\calR}{\mathcal{R}}

\newcommand{\one}{\mathbf{1}}
\newcommand{\zero}{\mathbf{0}}

\newcommand{\C}{\mathbb{C}}
\newcommand{\R}{\mathbb{R}}
\newcommand{\Hil}{\mathbb{H}}

\newcommand{\dimension}{\textrm{dim}}
\newcommand{\TL}{\textrm{TL}}
\newcommand{\tr}{\textrm{tr}}
\newcommand{\QL}{\textrm{\bf{QL}}}

\DeclareMathOperator{\Id}{Id}

\numberwithin{equation}{section}

\newtheorem{theorem}[equation]{Theorem}

\newtheorem{proposition}[equation]{Proposition}

\newtheorem{definition}[equation]{Definition}

\newtheorem{thm}{Theorem}
\newtheorem{lemma}[thm]{Lemma}
\newtheorem{defn}[thm]{Definition}
\newtheorem{rmk}[thm]{Remark}
\newtheorem{coro}[thm]{Corollary}
\newtheorem{examp}[thm]{Example}
\newtheorem{prop}[thm]{Proposition}

\newcommand{\bt}{\begin{thm}}
\newcommand{\et}{\end{thm}}
\newcommand{\bl}{\begin{lemma}}
\newcommand{\el}{\end{lemma}}
\newcommand{\br}{\begin{rmk}}
\newcommand{\er}{\end{rmk}}
\newcommand{\bc}{\begin{coro}}
\newcommand{\ec}{\end{coro}}
\newcommand{\bd}{\begin{defn}}
\newcommand{\ed}{\end{defn}}
\newcommand{\bex}{\begin{examp}}
\newcommand{\eex}{\end{examp}}
\newcommand{\beq}{\begin{equation}}
\newcommand{\eeq}{\end{equation}}
\newcommand{\bp}{\begin{prop}}
\newcommand{\ep}{\end{prop}}

\newcommand{\pair}[1]{\langle #1 \rangle}

\begin{document}
\title[Third Life]
{The Third Life of Quantum Logic: \\ Quantum Logic Inspired by Quantum Computing}

\author{J. Michael Dunn}
\email{dunn@indiana.edu}
\address{School of Informatics and Computing\\
    Indiana University\\
    Bloomington, IN 47405\\
    U.S.A.}

\author{Lawrence S. Moss}
\email{lsm@cs.indiana.edu}
\address{Department of Mathematics\\
    Indiana University\\
    Bloomington, IN 47405\\
    U.S.A.}

\author{Zhenghan Wang}
\email{zhenghwa@microsoft.com}
\address{Microsoft Station Q\\Elings Hall 2237\\
    University of California\\
    Santa Barbara, CA 93106\\
    U.S.A.}

\begin{abstract}
We begin by discussing  the history of quantum logic,  dividing it into three eras or ``lives.''
The first life has to do with Birkhoff and von Neumann's algebraic approach in the 1930's.  The second life has to do with attempt to understand quantum logic as logic that began in the late 1950's and blossomed in the 1970's.  And the third life has to do with recent developments in quantum logic coming from its connections to quantum computation.  We discuss our own work connecting quantum logic to quantum computation (viewing quantum logic as the logic of quantum registers storing qubits), and make some speculations about mathematics based on quantum principles.   
\end{abstract}

\maketitle

\section{History}

Modern classical logic began with Boole (1847), who had two interpretations of
the elements in his algebra of logic. The first interpretation was that they
were classes; the second was that they were propositions.  He
connected the two, saying that for purposes of inference a proposition could
be regarded as a class.\footnote{Boole (1847) spoke somewhat abstractly of
classes of \textquotedblleft conceivable cases and conjunctures of
circumstances,\textquotedblright\ whereas Boole (1854) took a more concrete
temporal interpretation, speaking of classes of \textquotedblleft
times\textquotedblright\ (calling these \textquotedblleft
durations\textquotedblright). See Kneale (1956).}

We see in Boole a prescient anticipation of the contemporary mathematization
of propositions as sets (of \textquotedblleft cases,\textquotedblright%
\ \textquotedblleft occasions,\textquotedblright\ \textquotedblleft
times,\textquotedblright\ \textquotedblleft possible worlds,\textquotedblright%
\ \textquotedblleft situations,\textquotedblright\ \textquotedblleft
set-ups,\textquotedblright\ \textquotedblleft states,\textquotedblright%
\ whatever), sometimes called  \textquotedblleft UCLA
propositions.\textquotedblright\ Conjunction is interpreted as intersection,
disjunction as union, and negation as complement (relative to a given
underlying set of possible ).  This way of looking at propositions can be generalized to include other non-classical logics, including quantum logic, though we will see that only conjunction remains in its original form.

\subsection{The first life of quantum logic: Birkhoff and von Neumann}

Quantum logic began with Birkhoff and von Neumann (1936) when they published
their pioneering paper titled \textquotedblleft The Logic of Quantum
Mechanics.\textquotedblright\ \footnote{That
paper was anticipated by von Neumann's 1932 book on the mathematical
foundations of quantum mechanics. \ There in section 5, chapter 3, he observed
that the projections defined on a Hilbert space could be regarded as
representing experimental propositions concerning the properties of a quantum
mechanical system. \ Projections correspond to closed subspaces.}

They point out that in classical dynamics, the state of a single particle can
be described as a sextuple $\pair{r_{1},\ldots,r_{6}}$ of real numbers, where the
first 3 components specify its position and the second 3 components specify
its momentum. $\ $The \textquotedblleft phase space\textquotedblright\ for
$n\ $particles can thus by thought of as the product set
$%
%TCIMACRO{\U{211d} }%
%BeginExpansion
\mathbb{R}
%EndExpansion
^{6n}$. Any subset of $%
%TCIMACRO{\U{211d} }%
%BeginExpansion
\mathbb{R}
%EndExpansion
^{6n}$ can be thought of as an \emph{event}, or proposition. And conjunction,
disjunction, and negation can be understood as Boole's operations on sets.

But anyone who knows anything about quantum mechanics has heard of the
Heisenberg Uncertainty principle, which says that one cannot simultaneously
determine both the position and momentum of a particle. \ Birkhoff and von
Neumann accordingly move to the more complicated phase-space on which they
build their quantum logic. There is a correspondence with classical dynamics
and Boole, except not every set of states determines a proposition -- only the
closed subspaces. The conjunction $\wedge$ of two subspaces is their set
intersection, but their disjunction $\vee$ is the closure of their
span.\footnote{This enlarges the union in two distinct ways. \ First by adding
all linear combinations (the \textquotedblleft span\textquotedblright), and
secondly by adding all limit points (the \textquotedblleft
closure\textquotedblright).} And the negation of a subspace is the set of
states that are \textquotedblleft orthogonal\textquotedblright\ ($\alpha
\perp\beta$) to every state in the subspace, where $\alpha\perp\beta$ means
that the \textquotedblleft inner product\textquotedblright\ $\alpha
\bullet\beta=0$.

\subsection{Boolean algebras and their generalizations.}
In this section we quickly review some algebraic structures that have naturally arisen in the study of classical logic and quantum logic.

A \textit{Boolean algebra} is a special kind of bounded distributive lattice
where every element $x$ has a complement $\sim\!x$.  Let us go through this a
step at a time. A \textit{lattice} can be defined as a partially ordered set
$(L,\leq)$ where for any $x$, $y$ $\in L$, there exists a greatest lower bound
$x\wedge y$ $\in L$ and a least upper bound $x\vee y$ $\in L.$ 
(There are several additional laws which we omit.)  
Think of $x,y$
as propositions, $\leq$ as entailment, $\wedge$ as conjunction, and $\vee$ as disjunction.

The lattice is \textit{bounded} if it has both a least element $0$ and a
greatest element $1$. It is \textit{complemented} if for every element $x$
there is an element $x^{\prime}$ such that $x\wedge x^{\prime}=0.$ Complements
are unique in a Boolean algebra, so we can introduce a unary operation $\sim$
that takes $x$ to its complement. It can then be shown that $\sim$ is
\textit{order inverting} (if $x\leq y$ then $\sim\!y\leq\,\sim\!x$) and of
\textit{period two} ($\sim\sim\!x=x$). In the context of a lattice these are
equivalent to the \textit{De Morgan Laws}: $\sim\!(x\wedge y)=$ $\sim
\!x\,\vee\sim\!y$ and $\!\sim(x\vee y)=$\ $\sim\!x\,\wedge\sim\!y.$ In a
Boolean algebra then we have $x\vee\sim\!x=1.$ A lattice is
\textit{distributive} if $x\wedge(y\vee z)\leq(x\wedge y)\vee(x\wedge z).$
(The converse is true in any lattice.)

A unary operation $\sim$  on an arbitrary lattice is an
\textit{orthocomplement} if it is of period two, order inverting, and
$\sim\!x$ is the complement of $x$.\ Orthocomplements are not necessarily
unique outside of the context of a distributive lattice. \ An
\textit{ortholattice} is a structure $(L,\leq,\sim)$\ where $(L,\leq)$ is a
lattice and $\sim$ is an orthocomplement.

An ortholattice is a generalization of a Boolean algebra in that it does not
need to be distributive.

Let us first consider three critical laws. \ There are various equivalent ways
to state the modular law, and we choose:%
\[
\text{(Modularity) if }z\leq x,\text{then }x\wedge(y\vee z)\leq(x\wedge y)\vee
z.
\]
The modular law holds in classical logic, and in fact it holds without any
conditions, since in the context of a lattice the consequent ($x\wedge(y\vee
z)\leq(x\wedge y)\vee z$) comes unconditionally from distribution (and is in
fact equivalent). Modularity can also be stated unconditionally as:%
\[
\text{(Unconditional Modularity) }x\wedge(y\vee\lbrack x\wedge z])\leq(x\wedge
y)\vee z.
\]
This is not strictly an equation but we can treat all weak inequalities as
equations in virtue of the general lattice equivalence $s\leq t$ iff $s\wedge
t=s$. This means that the class of modular lattices is equationally definable.

Birkhoff and von Neumann noted that the distributive law fails in their
quantum logic, but somewhat surprisingly they also note that the modular law
also fails. \ Instead there is a weaker \textquotedblleft orthomodular
law\textquotedblright:%

\[
\text{(Orthomodularity) if }z\leq x,\text{then }x\wedge(\sim\!x\vee z)\leq
z.\footnote{This has the philosophically memorable equivalent: $x\wedge
(\sim\!x\vee(x\wedge y))\leq y$, \linebreak which has led to regarding
$\sim\!x\vee(x\wedge y)$ as a conditional -- the so-called \textquotedblleft
Sasakai hook,\textquotedblright\ named after its discoverer.}
%\par
%{}}%
\]

An \textit{orthomodular lattice} is an ortholattice in which the orthomodular
law holds. \ This leads to a linguistically confusing but important
distinction between an orthomodular lattice and a modular ortholattice.
Modular ortholattices are special kinds of orthomodular lattices. It is
interesting that Birkhoff and von Neumann (1936) took the former and not the
latter as part of their logic of quantum mechanics. R\'{e}dei points out that
the they prefer the modular law because of its fit with a generalization of
classical probability theory.\footnote{R\'{e}dei (2007) contains an
interesting discussion of this, as do other publications by R\'{e}dei. \ See
particularly R\'{e}dei (2005) regarding the background correspondence from von
Neumann to Birkhoff. \ }

Our reason for liking the modular law was different and had to do with wanting
a generalization of the register of bits in a classical computer, so as to
have a quantum register of qubits. \ We focused on finite registers and it is
well-known that finite dimensional orthomodular lattices are modular.

There is an absraction under which one can fit both probability and
dimension.  A standard (Kolmgorov) requirement on  a probability function

\[
p(a\vee b)=\ p(a)+p(b)-p(a\wedge b)
\]
can trivially be restated and generalized (putting a general function $f$ for
$p$) as
\[
d(a)+d(b)=d(a\vee b)+d(a\wedge b).
\]
Birkhoff (1940) calls such a real-valued function a \textquotedblleft
valuation\textquotedblright\ and shows that the existence of a strictly
monotonic valuation on a lattice implies implies that the lattice is modular.
Birkhoff observes that both probability and dimension are valuations, and both
are monotonic. \ Dimension is obviously also strictly monotonic $(a<b$ implies
$d(a)<d(b))$, and so it seems is probability when it is taken in its logical
interpretation.\footnote{Von Neumann seems to have gone back and forth on how
he interpreted probability (frequency or logical), but about this time seemed
to favor logical probability. \ See Redei (2005).}

\subsection{The second life of quantum logic: quantum logic as logic.}
From the late 1950's, and especially in the 1970's and 80's, quantum logic had
a second life. \ As is said by Dalla Chiara and Giuntini after discussing
Birkhoff and von Neumann (1937): \textquotedblleft Only twenty years later,
after the appearance of George Mackey's book \textit{Mathematical Foundations
of Quantum Theory} [Mackey, 1957], one has witnessed a `renaissance period'
for the logico-algebraic approach to QT. This has been mainly stimulated by
the contributions of Jauch, Piron, Varadarajan, Suppes, Finkelstein, Foulis,
Randall, Geechie, Gudder, Beltrametti, Cassinelli, Mittelstaedt and many
others.\textquotedblright\ The main topic of interest regarding quantum logic
regarded the novelty of yet another non-classical logic, and how it compared
with intuitionistic logic (the main alternative to non-classical logic at that
point). There was also relatively great interest (compared to Birkhoff and von
Neumann, and now) about how it should best be conceived (orthomodular lattice,
many-valued logic, etc.), and following this in the standard logical issues of semantics,
proof-theory, completeness, and decidability. There was a strong preference
for the orthomodular approach, and that gave a proof-theory and a semantics
(the closed subspaces of a Hilbert space), but connecting the two has proved
impossible. It turns out that the lattice of closed subspaces of a Hilbert
space satisfies additional laws, even those that can be stated as equations
such as the \textquotedblleft Orthoarguesian law\textquotedblright%
\footnote{This apparently was an unpublished idea of Alan Day, and the proof
was first presented by Greechie (1983). \ See Dalla Chiara, Giuntini, and
Greechie (2004)}.    Another problem is that the axiom system can be given in
so-called \textquotedblleft Hilbert-style\textquotedblright\ by translating
the axioms for orthomodular lattices into a more standard logical
formalism,\ but to our knowledge no one has yet succeeded in giving an equivalent
cut-free Gentzen version (which many people think is the gold-standard
approach to proof-theory) of orthomodular logic or modular
orthologic.\footnote{See Nishimura (2009) for presentation and history of
cut-free Gentzen systems for \textquotedblleft minimal quantum
logic\textquotedblright\ (what we are calling orthologic) and its history.
\ See also Egly and Tompits (1999). Chiara and Giuntini (2002) in sec. 17 (by
G. Battilotti and C Faggian) discuss a Gentzen system for orthologic developed
by Sambin, Battilotti, and Faggian that has a cut-free formulation, but they
do not address orthomodular logic or modular orthologic.}

\subsection{The third life of quantum logic: quantum logic based on quantum computation.}
Dalla Chiara and Giuntini speak of a \textquotedblleft
Renaissance,\textquotedblright\ which of course literally means
\textquotedblleft rebirth,\textquotedblright\ i.e., a second life.  We believe
that quantum logic now has a \textquotedblleft third life,\textquotedblright%
\ inspired by quantum computing. And subtly different algebraic structures
arise (with some of the same open questions, but new chances at solving them).

Although Richard Feynman had first suggested the idea of a quantum computer to
simulate quantum processes faster than might be done by a classical computer,
it was not until 1985 that David Deutsch published a paper describing a
general purpose (universal) quantum computer. Deutsch modified the classical
Turing machines to make equivalents among other differences put qubits in
place of the standard binary digits that appear on the tape of a Turing
machine. The equivalent use of quantum gates has become the much more
customary way of characterizing quantum computation.\ This work was once
purely theoretical, but after 1994, when Peter Shor gave his famous algorithm
for efficiently factoring numbers into prime, the idea of a quantum computer
began to take on real practical significance. This is because of the widely
used RSA encryption scheme that depends on the difficulty of factoring large
numbers into their prime components.

From its early years logic has been linked to computation. Leibniz's great
achievement was to combine the idea of a \textquotedblleft lingua
universalis\textquotedblright\ with a \textquotedblleft calculus
raciocinator\textquotedblright. The two together facilitate \textquotedblleft
blind thinking,\textquotedblright\ as Leibniz termed it, since reasoning is
reduced to arithmetic calculation.\textquotedblright The link between
classical computing and classical logic is often taken for granted.\ Many
standard classical textbooks contain both, e.g. Kleene's (1950)
\textit{Introduction to Metamathematics.} But the use of classical logic to
describe and design circuits is not even mentioned, whereas this has become
almost the standard approach to thinking about quantum logics.

In a classical computer, data is stored as a \textquotedblleft
string\textquotedblright\ of bits in a register. Registers come in various
sizes, thus a 64-bit register contains strings of length 64. The
\textquotedblleft register space\textquotedblright\ can be viewed as the
direct product of the 2-element Boolean algebra, i.e., the set of $n$-element
sequences of 0s and 1s. This can be viewed as a Boolean algebra itself by the
direct product construction, defining $\wedge$,$\vee,-$ component wise, e.g.,
$-\langle b_{1},\ldots,b_{n}\rangle=\langle-b_{1},\ldots,-b_{n}\rangle.$

What is the logic of the classical $n$-bit register? Is it the same as the
logic of the $1$-bit register (classical logic) or not? This is answerable in
two steps. \ We first take classical propositional calculus, and form its
\textquotedblleft Lindenbaum Algebra\textquotedblright\ by the ``Method of
Abstraction." \ We thus put two provably equivalent formulas $\varphi$ and
$\varphi^{\prime}$ into the same equivalence class $[\varphi]=[\varphi
^{\prime}],$ and we then define operations on these equivalence classes using
the sentential operators, e.g., $-[\varphi]=[-\varphi].$ \ For classical
propositional calculus this gives a Boolean algebra in which the equivalence
class of the theorems $=1$. The second step is to invoke a form of the
Representation Theorem for Boolean algebras (Stone 1935): Every Boolean
algebra is isomorphic to a subdirect product of the $2$-element Boolean
algebra ($1$-bit register). Combining these ideas, classical propositional
logic can easily be shown to be the logic of the $n$-bit register (not just
the $2$-bit one).

The qubit is a \textquotedblleft quantum bit\textquotedblright. Unlike the
classical bit, $0$ and $1$ are just two of infinitely many possible states of
the qubit. The state of a qubit is the \textquotedblleft
superposition\textquotedblright\ (linear combination) $\alpha 0\rangle
+\beta 1\rangle$ (where $\alpha,\beta$ are complex numbers representing
\textquotedblleft amplitudes\textquotedblright--amplitudes squared give
probabilities).\footnote{This is usually written in the Dirac notation as
$\alpha|0\rangle+\beta|1\rangle$, but we will not be so fussy in our
motivating explanations here.} The state of a qubit can be described as a
vector $(\alpha,\beta)$ in the two-dimensional complex vector space $%
%TCIMACRO{\U{2102} }%
%BeginExpansion
\mathbb{C}
%EndExpansion
^{2}$. The special states $0$ and $1$ are known as the computational basis
states, and form an orthonormal basis for this vector space. According to
quantum theory, when we try to measure the qubit in this basis in order to
determine its state, we get either $0$ with probability $|\alpha|^{2}$ or $1$
with probability $|\beta|^{2}$. This motivates requiring that $|\alpha|^{2}%
+|\beta|^{2}=1$. (This is a Probability Sum Rule for disjoint events).

Quantum registers contain qubits (quantum bits). The $2$-dimensional space of
the complex numbers $%
%TCIMACRO{\U{2102} }%
%BeginExpansion
\mathbb{C}
%EndExpansion
^{2}=%
%TCIMACRO{\U{2102} }%
%BeginExpansion
\mathbb{C}
%EndExpansion
\oplus%
%TCIMACRO{\U{2102} }%
%BeginExpansion
\mathbb{C}
%EndExpansion
$ can be thought as a quantum register containing a single qubit, and all the
pairs of complex numbers in it are then thought of as states of that qubit.
The $n$-qubit register
$%
%TCIMACRO{\U{2102} }%
%BeginExpansion
\mathbb{C}
%EndExpansion
^{2^{n}}$ can then be defined inductively as $%
%TCIMACRO{\U{2102} }%
%BeginExpansion
\mathbb{C}
%EndExpansion
^{2}\otimes%
%TCIMACRO{\U{2102} }%
%BeginExpansion
\mathbb{C}
%EndExpansion
^{2^{n-1}}$, i.e., $%
%TCIMACRO{\U{2102} }%
%BeginExpansion
\mathbb{C}
%EndExpansion
^{2}\otimes%
%TCIMACRO{\U{2102} }%
%BeginExpansion
\mathbb{C}
%EndExpansion
^{2}\otimes\cdots\otimes%
%TCIMACRO{\U{2102} }%
%BeginExpansion
\mathbb{C}
%EndExpansion
^{2}$ ($n$-times). \ It turns out that unlike the analogous case with
classical logic, the logic of the $n$-qubit register is generally different
that the logic of the $1$-bit quantum register, and indeed the logic of the
$n$-qubit register is always different from the logic of the $m$-bit quantum
register when $m\neq n.$ This was shown in Dunn, Hagge, Moss, and Wang (2005), and the result
was improved by Hagge (2007) who showed for all $m\neq n$, the logics of 
$\mathbb{C}^m$ and $\mathbb{C}^n$ differ.
(Note that the superscript here is $n$ and not $2^{n}$.) This leads us
to wonder whether every  subdirectly irreducible modular lattice is
isomorphic to the lattice of subspaces of some $%
%TCIMACRO{\U{2102} }%
%BeginExpansion
\mathbb{C}
%EndExpansion
^{n}$.
If it were, then by   Birkhoff's Subdirect Product Theorem,
every modular ortholattice would be isomorphic to a
subdirect product of such lattices.

We titled the conference \textquotedblleft Quantum Logic Inspired by Quantum
Computing" (QLIQC, pronounced \textquotedblleft click\textquotedblright), but
it turns out it might just as well have been \textquotedblleft Quantum Logic
Inspired by Quantum Categories" in terms of the talks given (and the
subsequent papers published in this volume). \

\section{Quantum logic inspired by quantum computing}

The BB$84$ private key protocol (Bennett and Brassard 1984), Shor's algorithm (Shor 1994), and Hastings' additivity counterexamples (Hastings 2009) are all 
pieces of evidence  that quantum information theory is strictly richer than classical information theory.
The  attempt to build 
 a useful quantum computer has begun and rekindled interest  in quantum mechanics at all levels: philosophical, mathematical and physical. Unlike the  construction of the classical computer, to build a quantum computer might require new physics such as non-abelian topological order (Freedman et al 2003).

Ever since its appearance, quantum mechanics presents great conceptual difficulty, even for the most brilliant minds.  In quantum mechanics, the wave function of a state is a complete description of the physical state,
 and the Schrodinger equation is a deterministic evolution of the state.  When the measuring apparatus is included into the system, the measurement of a quantum system is a deterministic process for the composite system with complete description.  Yet our best interpretation for the measurement result is still probabilistic.
Probability is usually related to insufficient knowledge.  The mismatch of a complete description of a quantum system with the probabilistic interpretation lies at the heart of the debate.  Maybe humans are innately not able to apprehend a quantum state.  But the emergence of numbers seems to suggest otherwise.  Children's counting ability is arguably primitive and dormant, and only through education is the number fully developed into counting with numbers.  Historically tally seems to come first, then counting, and finally abstract numbers.  An important step in the emergence of numbers is the separation of things to be counted from their associated symbols.  Quantum information is taking this step right now.  The qubit is the abstraction of $2$-level quantum systems, therefore it is not an electron spin;
similarly, the number one is not an apple.  Hence the qubit likely will play an important role in the evolution of numbers.

\subsection{Quantum logic of qubits}

Qubits are the currency for quantum computing.  Their states are represented by non-zero vectors of the Hilbert space ${(\mathbb C^2})^{\otimes n}$.  In this section, we will examine the quantum logic of qubits following the ideas of G.~ Birkhoff and J. ~von Neumann (1936).

\subsubsection{Quantum logic determines dimension}

Given a Hilbert space $\mathbb{H}$, let $L_c(\Hil)$ be the lattice of closed subspaces of $\Hil$.  Closed subspaces are quantum events, so they are quantum analogues of propositions.  We will use $\zero,\one$ to denote the $0$-subspace and $\Hil$, respectively.  The meet $\wedge$ of two subspaces is the set intersection, and the join $\vee$ the closure of their span.  For any closed subspace $p$, its negation 
${p}^{\bot}
$ is the orthogonal complement.  It is well-known that $L_c(\Hil)$ is an orthomodular lattice and modular if and only if $\Hil$ is finite dimensional.  Propositional formulas consist of alphabet symbols, parenthesis, and connectives $\wedge, \vee$ and $\;\bar{}\;$.  Let $\{u_i\}$ be a collection of alphabet symbols, and $\{p_i\}$ be a collection of closed subspaces.  Given a well-formed formula (wff) $\phi(u_i)$, the evaluation $\phi(p_i)$ is the subspace resulting from substituting each $p_i$ into $u_i$ and performing the corresponding 
operations.  A wff $\phi(u_i)$ is a
\emph{tautology} of $L_c(\mathbb H)$ if for all evaluations $\phi(u_i)=\one$.  We will also call any equation of terms $s = t$ in which for all evaluations $s = t$ a \emph{tautology}.

\begin{definition}
Given a Hilbert space $\Hil$, the quantum logic ${\bf QL}(\Hil)$ is
the set of all tautologies of $L_c(\mathbb H)$.
\end{definition}

\begin{theorem}
Quantum logic ${\bf QL}(\Hil)$ determines the dimension of $\Hil$.
\end{theorem}

Note that the modular law separates infinite dimensional Hilbert spaces from finite dimensional ones.  Then the dimension of a finite dimensional Hilbert space is determined by its quantum logic (Dunn, Hagge, Moss, and Wang 2005, Hagge 2007).

For notational ease, we will denote ${\bf QL}(({\mathbb C^2})^{\otimes n})$ by ${\bf QL}(2^n), n=0,1,\cdots$.  To understand the differences between these logics better, we will exhibit tautologies that distinguish them.  For  $n=0$, the quantum logic ${\bf QL}(1)$ is just the classical propositional logic.  The distributive law holds in ${\bf QL}(1)$, but fails in any ${\bf QL}(2^n), n\geq 1$.  Therefore, the distributive law is a salient feature of classical logic.  We will explore the failure of distributive law in ${\bf QL}(2^n), n\geq 1$ systematically to arrive at increasingly weakened tautologies.  The first such tautology was the $m$-distributive law:

$$x\wedge (\vee_{i=0}^m y_i)=\vee_{i=0}^m (x\wedge(\vee_{j\neq i} y_j)).$$

It is proven (Huhn 1972) that the $m$-distributive law holds if and only if $\dimension (\Hil)\leq m$.  Dunn, Hagge, Moss, and Wang (2005), and Hagge (2007) found another sequence of such tautologies.  For simplicity, we will consider only the qubits here.  As a bonus of our new tautologies, we will see that ${\bf QL}({2^n})$ has no finite universal test sets when $n\geq 1$.

Two closed subspaces $a,b$ are equal if and only if $(a\vee b)\wedge (\bar{a}\vee \bar{b})=\zero.$  To see this, if $a=b$, obviously $(a\vee b)\wedge (\bar{a}\vee \bar{b})=\zero$.  If $a\neq b$, then either $a\wedge b\neq a$ or $a\wedge b\neq b$.  Without loss of generality, we assume $a\wedge b\neq a$.  Then the complement of $a\wedge b$ in $a$, denoted as $\overline{a\wedge b}^a$, is not $\zero$.  But $\overline{a\wedge b}^a \subset \overline{a\wedge b}=\bar{a}\vee \bar{b}$ and $\overline{a\wedge b}^a\subset a\vee b$.  Hence $(a\vee b)\wedge (\bar{a}\vee \bar{b})\supset \overline{a\wedge b}^a \neq \zero.$

Given three subspaces $p,q,r$, let
$a=p\vee(q\wedge r)$ and $b=(p\vee q)\wedge (p\vee r)$, and then
define
$$\alpha(p,q,r)=(a\vee b)\wedge (\bar{a}\vee \bar{b}).$$
Note that $a\leq b$, it follows that
$\alpha(p,q,r)=b\wedge \bar{a}=[(p\vee q)\wedge (p\vee r)]\wedge
[\bar{p}\wedge(\bar{q}\vee \bar{r})]\subseteq \bar{p}.$
The distributive law holds if and only if
$\alpha$ is always ${\bf 0}$.  Therefore,
if $\alpha$ does not vanish for some choice of $p,q,r$ in
a Hilbert space $\Hil$, then the distributive law is not in ${\bf QL}(\Hil)$.  Therefore, we will call $\alpha(p,q,r)$ the {\it distribution test formula}.

From $\alpha(p,q,r)\subset \bar{p}$, we deduce $\dimension(\alpha(p,q,r))\leq \dimension(\Hil)-\dimension(p)$.  In Dunn, Hagge, Moss, and Wang (2005), a direct computation shows $\dimension(\alpha(p,q,r))\leq \dimension(p)$.  Hence $\dimension(\alpha(p,q,r))\leq \frac{\dimension(\Hil)}{2}.$

To define our tautology, we define the restriction of a wff $\phi(u_i)$ to a term $\beta$, denoted by
$\phi|_{\beta}$:  first using the De Morgan
law, we assume that all negations \;$\bar{}$\; are applied to
single variables.   Next,  each variable $u_i$ and  its complement
$\bar{u_i}$ are replaced
by $u_i\wedge \beta$ and $\overline{u_i\wedge \beta}\wedge \beta$, respectively.
Inductively, we define
$$\alpha^m(p_m,q_m,r_m)=\alpha|_{\alpha^{m-1}}({p_m},{q_m},{r_m}),$$
and $\alpha^1(p_1,q_1,r_1)=\alpha(p_1,q_1,r_1),
\alpha^{m-1}=\alpha^{m-1}(p_{m-1},q_{m-1},r_{m-1}).$
Therefore,
$$\textrm{dim}(\alpha^m(p_m,q_m,r_m))
\leq \frac{\textrm{dim}(\alpha^{m-1}(p_{m-1},q_{m-1},r_{m-1}))}{2}
\leq \cdots \leq \frac{\dimension(\Hil)}{2^m}.$$

In ${\bf QL}(2^n)$,
$\textrm{dim}(\alpha^{n+1})\leq \frac{2^n}{2^{n+1}}<1$, so
$\alpha^{n+1}={\bf 0}$ which gives a tautology in $\QL({2^n})$, which is also true for any $i\leq
n$.  To show it is not true for ${\mathbb C}^{2^{n+1}}$, we notice
that if $p, q, r$ are different subspaces of dimension
$\frac{m}{2}$ of $\C^m$ and each pair has trivial intersection in ${\mathbb C}^m$, then
$\textrm{dim}(\alpha(p,q,r))=\frac{m}{2}$ if $m$ is even.  By choosing
subspaces in ${\mathbb C}^{2^{n+1}}$ this way, we have
$\textrm{dim}(\alpha^{n+1})=\frac{2^{n+1}}{2^{n+1}}=1$.

\begin{definition}
A set of closed subspaces in
${\mathbb C}^m$ is called \emph{a universal test set} for $\QL(\C^m)$ if
the truth of any tautology is
determined by the evaluations of the subspaces in this set.
\end{definition}

\begin{proposition}
There are no finite universal test sets for ${\bf QL}({\mathbb
C}^m), m\geq 2$.
\end{proposition}

To see this, consider
the distribution testing formula $\alpha(p,q,r)$. For simplicity, we
will only give the details for $m=2$.  In order
for the distribution testing formula $\alpha(p,q,r)$ to
fail, $p, q, r$ must be three distinct lines.
In order for $\alpha(\alpha(\alpha(\alpha(p,q,r),p,s),q,s),r,s)$
to fail, $p, q,r,s$ must
be distinct lines.
Continuing in this manner, we can build a complicated formula
$\gamma$, the failure of which means that the $k$ subspaces
$p,q,\cdots$ are distinct lines.  Since $k$ is arbitrary,
no finite set of lines will
falsify every invalid formula.  This argument works for
any ${\mathbb C}^m, m\geq 2$.

For each $n$-qubit, we have found two tautologies which are not in any qubits $m$ such that $m<n$: the $2^n$-distributive law and the iterated
distribution test formula.  If each law is added to the modular lattice axioms, are the resulting 
axioms sets equivalent? We leave this as an open problem.

\subsubsection{Decidability}

Quantum logic for general modular ortholattice is undecidable.  Dunn, Hagge, Moss, and Wang (2005) observed that the quantum logic of a finite dimensional
 Hilbert space is decidable. The decidability ${\bf QL}({\C}^m)$ is reduced to the decidability of $\R$.  The idea is to associate a matrix $M_p$ to each subspace $p$ so that the kernel of $M_p$ is $p$.  Then new matrix variables are introduced to construct a formula $M_{\phi}$ so that a wff $\phi$ is a tautology if and only if $M_{\phi}=0$.  This procedure is illustrated for the join in Dunn, Hagge, Moss, and Wang (2005).  The easy cases of meet and negation can be done as follows:
$$r=p\wedge q, \forall u(M_p u=0\wedge M_q u=0 \Leftrightarrow M_r u=0).$$
$$p=\bar{q}, \forall v (M_p v=0\Leftrightarrow (\forall u (M_q u=0\Rightarrow <v,u>=0))).$$
Taking all these observations together, we conclude:

 \begin{theorem}
 The first-order theories of ${\bf QL}({\C}^m)$ are uniformly decidable.
 \end{theorem}

Since decidability of ${\bf QL}({\C}^m)$ implies its axiomatizability, can ${\bf QL}({\C}^m)$ be axiomatized
with finitely many schemas?  Is it sound? Is it complete?  Are modular ortholattice axioms plus $n$-distributivity or the iterated distribution test formula sufficient to axiomatize ${\bf QL}({\C}^m)$?  We believe that these are all interesting open problems.

We might also speculate on a connection between quantum logic and quantum computational complexity.  For example, if we choose a finite collection of subspaces of $\C^2$ including $\zero, \one$ that generate a sublattice, then what is the computational complexity for the satisfiability?  In particular, if we add one $p$ whose normalized dimension is $\frac{1}{2}$ to $\zero, \one$, does the computational complexity depend on the choice of $p$?  Does quantum computer have any advantage over classical computers for those problems?

%Another direction will be the relation to non-classical logics such as multi-valued logics.

\subsection{Qubit continuous geometry}

Birkhoff and von Neumann proposed continuous geometry as quantum propositional logic.  In this section, we will focus on a particular continuous geometry---qubit continuous geometry.  This turns out to be the famous type $\Pi_1$ hyperfinite factor $\calR$ in von Neumann algebra theory.  Through the study of type $\Pi_1$ factors, V.~ Jones discovered his famous representation of the braid groups and polynomial invariants of knots.  Jones' representation of braid groups are used to describe new particle statistics and
are therefore playing a pivotal role in the topological approach to quantum computing.

Let $V,W$ be two Hilbert spaces.  Note that neither $V$ nor $W$ is canonically a subspace of $V\otimes W$.  But the lattices $L_c(V)$ and $L_c(W)$ are canonical sublattices of $L_c(V\otimes W)$ by including $p\subseteq V$ or $W$ into $V\otimes W$ as $p\otimes W$ or $V\otimes p$, respectively.  It follows that $\bf{QL}(V)$ and $\bf{QL}(W)$ are canonically subsets
 of $\bf{QL}(V\otimes W)$. Therefore, quantum logics of qubits form a compatible decreasing sequence:
$$ {\bf QL}(1)\supset {\bf QL}(2)\supset
{\bf QL}(4)\supset\cdots \supset
{\bf QL}({2^n})\supset {\bf QL}({2^{n+1}})\supset\cdots.$$
 How to describe their intersection ${\bf QL}({\infty})$?  As remarked in Dunn, Hagge, Moss, and Wang (2005), the intersection ${\bf QL}({\infty})$ is not the quantum logic of any infinite dimensional Hilbert space because it contains the modular law.

\subsubsection{Limit of $\bf{QL}(2^n)$}

The normalized dimension of a subspace $p\subseteq V$ is $d_V(p)=\frac{\dimension(p)}{\dimension(V)}: L_c(V)\rightarrow [0,1].$
The lattice $L_c(2^n)$ with the normalized dimension $$d_{{(\C^2)}^{\otimes n}}: L_c({(\C^2)}^{\otimes n})\rightarrow [0,1]$$
 is a metric space compatible with the inclusion $L_c(2^n)\subset L_c(2^{n+1})$.   %\marginpar{I changed the notation here.}
  Let $L_c(\infty)$ be their direct limit.  The ranges of dimensions are all rational numbers with power $2$ denominators.  Let $CG$ be its metric completion, then $CG$ is a continuous geometry: an irreducible complemented continuous modular lattice.   A continuous geometry is a projective geometry whose dimensions cover the unit interval $[0,1]$.

To relate this continuous geometry to the hyperfinite $\Pi_1$ factor $\calR$, we consider the sequence of matrix algebras:
$$ M_2(\C)\subset M_2(\C)\otimes M_2(\C)\subset \cdots  $$ with inclusion given by $A\rightarrow A\otimes \Id$.   The $*$-algebra limit is the hyperfinite $\Pi_1$ factor $\calR$.  Let $L_p(\calR)$ be the set of projectors in $\calR$: $p=p^{\dagger}, p^2=p$.  Using the identification of a subspace with a projection, we see that $L_p(\calR)=CG$.  The factor $\calR$ can be realized as a subalgebra of the bounded operators of a Hilbert space $\Hil$.  With this realization, a projector can be identified with the closed subspace $p\Hil\subset \Hil$---invariant vectors of $p$ in $\Hil$.  We define the partial order, meet and join of two projectors $p,q$ by $p\leq q$ if and only if $p\Hil \subseteq q\Hil$, $p\wedge q=$orthogonal projection onto $p\Hil \wedge q\Hil$, $p\vee q=$orthogonal projection onto $p\Hil \vee q\Hil$.  The negation of a projector $p$ is $\bar{p}=1-p$.  Let $\QL(\calR)$ be the tautologies over $L_p(\calR)$ or equivalently over the sublattice of $L_c(\Hil)$ consisting of invariant subspaces of a projector in $L_p(\calR)$.

\begin{theorem}
\begin{enumerate}
\item $\QL(\infty)=\QL(\calR)$.
\item $\QL(\calR)$ is decidable.
\end{enumerate}
\end{theorem}

It is shown by J.~ Harding that $\QL(CG)=\QL({\infty})$.  Since $L_p(\calR)=CG$, therefore $\QL(\calR)=\QL(\infty)$.  It follows that $L_p(\calR)$ is a modular lattice.  As also proved in this issue by J.~ Harding, $\QL(\calR)$ is decidable.  Therefore,  the intersection ${\bf QL}({\infty})$ is decidable, positively answering a  question in Dunn, Hagge, Moss, and Wang (2005).
 The lattice of projectors $L_p(\calR)$ is a natural generalization of qubit quantum logic agreeing with Hankel's principle of the preservation of formal laws.  The decidability of quantum logics of general continuous geometries seems to be open.

%A lattice is atomic if for every $a=\vee_i a_i\in L$ such that if $b< a_i$ for each $i$, then $b=0$.  Since a type $\Pi_1$ factor has no finite %projection except the trivial one, hence it is nonatomic.

\subsubsection{Temperley-Lieb algebra and Jones-Wenzl projector}

It is clear from the last subsection how to obtain a projector in $L_p(\calR)$ with its normalized dimension to be any rational number in the interval $[0,1]$.  In this section, we construct projectors with irrational algebraic normalized dimensions.  It is hard to imagine projectors with non-computable normalized dimensions.  For example, let $\omega$ be a Chaitin number in $[0,1]$ which encodes the halting problem for Turing machines.   Do the projectors with normalized dimensions $\omega$ have any relevance to reality?

To construct such projectors, we introduce the Temperley-Lieb (TL) algebras.  The TL algebra $\TL_n(A)$ at  $A=\pm i e^{\pm \frac{2\pi i}{4r}}, r\geq 3$ is a unital algebra with generators $1,e_1,e_2,\cdots, e_{n-1}$ and relations:
\begin{equation}\label{farcommunitivityTL}
e_i e_j=e_je_i,\;\;\; \textrm{if} \;\;\; |i-j| \geq 2,
\end{equation}
\begin{equation}\label{braidrelationTL}
e_i e_{i\pm1}e_i=\frac{1}{d^2}e_i,
\end{equation}
and
\begin{equation}\label{heckerelationTL}
e_i^2=e_i,
\end{equation}
where $d=-A^2-A^{-2}$.

$\TL_n(A)$ is also defined by the same presentation when $A$ is a variable.  In this case, they are matrix algebras over a function field and contain  some magic projectors, called Jones-Wenzl projectors:  each ${\TL}_{n}(A)$ contains a unique element $p_n$ characterized by:
$p_n^2 = p_n\neq 0$ and $e_i p_n = p_n e_i=0$ for all $1\leq i\leq
n-1$. Furthermore $p_n$ can be written as $p_n = 1 +U$ where
$U=\sum c_j h_j$, $h_j$ a product of $e_i$'s, $1 \leq i \leq
n-1$ and $c_j\in \C$.

$\TL_n(A)$ can be naturally included into $\TL_{n+1}(A)$, hence $p_i, i\leq n$ can be considered as elements of $\TL_{n+1}(A)$.  When $A$ is chosen as the complex numbers above,
 the Jones-Wenzl projectors are defined consecutively only for $n=1,\cdots, r-1$.  Moreover,
the TL algebra is not a matrix algebra.  For a fixed $r$, their quotients by $p_{r-1}$ considered as an element in each $\TL_n(A), n\geq r-1$ are matrix algebras.  We will call those matrix algebras, denoted by $J_n(A)$,  the 
\emph{Jones algebras}.  
The matrix summands of the decomposition of $J_n(A)$ are indexed by natural numbers $m=n\; mod \;2$.  Define the $n^{{th}}$ Chebyshev polynomial $\Delta_n (x)$
inductively by $\Delta_0 =1, \Delta_1 =x$, and $\Delta_{n+1} (x) = x \Delta_n(x) - \Delta_{n-1}(x)$.   Then the Markov trace on $J_n(A)$ is the weighted matrix trace $\tr^{Mark}(u)=\sum_m \Delta_m \tr(u)$, where $\tr(u)$ is the usual matrix trace for $u$.  The Jones algebra ${J}_{n}(A)$ is included into ${J}_{n+1}(A)$ naturally.  The limit of them is the hyperfinite $\Pi_1$ factor $\calR$.  The Markov trace is the limit of the normalized dimensions.  In the Jones algebras $J_n(A)$, the TL elements $e_i$'s are Hermitian, i.e. $e_i^{\dagger}=e_i$, hence $e_i\in L_p(\calR)$.  The same is true for each Jones-Wenzl projector $p_i,i=0,1,\cdots, p_{r-2}$.  Therefore, $e_i$'s and $p_j$'s are projectors.  The TL relations tell us that the images of $e_i$ and $e_j$ are orthogonal if $|i-j|\geq 2$, and the angle between the $i^{th}$ and $(i+1)^{th}$ is determined by $d$.  Their trace are given by
 $\tr(e_i)=\frac{1}{d^2}=\frac{1}{4}sec^2(\frac{\pi}{r}),i=1,2,\cdots,$ and $\tr(p_j)=\frac{\Delta_j(d)}{d^j},j=1,2,\cdots, r-2$ for any $r\geq 3$.  Hence the subspaces $e_i\Hil$ and $p_j\Hil$ have normalized dimensions $\{\frac{1}{4}sec^2(\frac{\pi}{r})\}$ and $\{\frac{\Delta_j(d)}{d^j}\}$ for any $r\geq 3$ and $j=1,2,\cdots, r-2$.

Do the projectors $\{p_j, e_i, j=1,2,\cdots, r-2, i=1,2,\cdots \}$ and the subspaces of normalized dimensions with power $2$ denominators form a universal test set for $\QL(\calR)$?

\subsection{Topological quantum computation}

Classical physics is the theoretical foundation for the construction of classical computers.  The failure of C.~ Babbage to complete his analytical engine in $1850$s was not due to some missing physics, but rather for engineering reasons.   The same might occur for some current proposals to build a quantum computer.  But one approach is different in this regard: topological quantum computation (Freedman et al 2003). The success of topological quantum computation hinges on the discovery of completely new particles: non-abelian anyons.  The defining feature of such hypothetical particles are their ground state degeneracy in the plane: suppose several non-abelian anyons are fixed in the plane, well-separated, their lowest energy states are still not unique.   There are different internal states of the system which cannot be determined by their positions and other local properties.

\subsubsection{Non-abelian Anyons}

The mathematical models of non-abelian anyons are unitary modular tensor categories, or the closely related unitary topological quantum field theories.  An
anyon is a simple object in the modeling unitary modular category (Wang 2010).   %\marginpar{Which category?}

The Jones algebras can be easily generalized to tensor categories, which are unitary modular tensor categories.  The Jones-Wenzl projectors are the simple objects of the resulting unitary modular tensor categories, hence anyons.  When $r\geq 4$, the projector $p_1$ is a non-abelian anyon.  Suppose there are $m$ of them in a plane at some fixed locations, well-separated, how do we describe their states?  Their Hilbert space decomposes into subspaces of different energies.  The lowest energy states are called the ground states.  They are not unique and form a Hilbert space of dimension exponential in $m$.  Therefore, we need exponential many states to describe the differences of $m$ non-abelian anyons.

In topological quantum computation, information is encoded into this vast degenerate groundstate manifold of non-abelian anyons and processed by braiding them around each other.  Anyons can be brought together to fuse, and the computational answer, encoded in the resulting anyon types, is an approximation of the Jones polynomial at $q=e^{\pm \frac{2\pi i}{r}}$.

 Jones theories are predicted to be realized in the fractional quantum Hall liquids.  For example, when $r=4$, the anyon $p_1$ is predicted to exist in the $\nu=5/2$ fractional quantum Hall liquid.  The ground states will be $2$-fold degenerate if a boundary condition is fixed for $4$ non-abelain anyons $p_1$.  Therefore, to \lq\lq count" the ground states of $4$ non-abelian anyons $p_1$ when $r=4$, two independent wave functions are required.

\subsubsection{Intrinsic entanglement}

The many-anyon state is a state with topological order.  Topological order is an internal, dynamical, non-local pattern of many-anyon systems  characterized by intrinsic entanglement.  In quantum mechanics, entanglement is defined with respect to a tensor decomposition of the relevant Hilbert space, which amounts to a measurement.  In a topological state, the ground state manifold has no natural tensor decompositions.  Therefore entanglement of topological states is intrinsic---a salient feature of topological order.

\section{Speculative Remarks on Quantum Analogues of Classical Objects}

The term \emph{quantum mathematics} is ambiguous.  One sense concerns the mathematics
needed to explain and work with quantum physics.
The other sense is more radical and has to do with an alternative to classical mathematics that is somehow founded on quantum principles.  We shall briefly discuss several such approaches.

\subsection{Quantum Cantor set}

Von Neumann algebra theory is an axiomatization of quantum mechanics and can be regarded as a non-commutative measure theory.  A von Neumann algebra $M$ with a normalized, normal trace $\rho$ is called a noncommutative  probability space.  An Hermitian operator $X$, a physical measurement, is a noncommutative random variable.  Its eigenvalues $Spec(X)=\{\lambda_i\} \subset \R$ are  in one-one correspondence with projections $P_{\lambda_i}$ to the eigenspaces.  The values $\rho(P_{\lambda_i})$  %\marginpar{Right?}
 of the projections $\{P_{\lambda_i}\}$ under $\rho$ define a probability distribution on $Spec(X)$.  Therefore, a type $\Pi_1$ factor is a natural noncommutative probability space.

In the classical world, logic, measure theory, and probability fit together via 
the Stone representation theorem.  For bit strings, their limit is the Cantor set. In this analogy, hyperfinite $II_1$ factor $\calR$ with the normalized trace is a probability theory for a quantum Cantor set.  What is a quantum Cantor set? What is a  quantum Boolean algebra? How to logicize type $\Pi_1$ factors?

\subsection{Quantum numbers}

%\marginpar{\bf find reference to Kripke}

Classical computers process bit strings, which can be regarded as numbers denoted by binary notation.  Numbers seem to
be rooted in a human's need to record the differences between say one sheep and a herd.  To understand the physical
properties of many quantum particle systems potentially
leads us to new numbers.  The states of quantum particles cannot be easily described by
numbers as they are given by wave functions.  Quantum computers are wave function processors.
Therefore, we argue that wave functions are quantum numbers.

A real number in base $2$ expansion can be considered as an array of bits $\{0,1\}$ on a bi-infinite Turing tape with a marked square for the separation of integral and fractional parts.  The squares are digit holders.  If the squares in the Turing tape correspond to basis elements of a Hilbert space, then a wave function can be thought as a generalization of numbers in two aspects:  the bases of a Hilbert space is not necessarily an array and the digits $\{0,1\}$ are replaced by any complex number.  There are axiomatizations of both natural numbers---Peano axioms and real numbers---Dedekind cuts.  Are there axiomatizations of wave functions---our proposed quantum numbers?
For caution, we mention a work of Dunn: if the first order Peano arithmetic is formulated with orthomodular quantum logic, then it has the same theorems as the first order Peano arithmetic (Dunn 1980).  Quantum mathematics is slippery business as Dunn showed that if one tries to formulate second order orthomodular quantum logic with a certain minimal principle of extensionality, one is doomed to failure in the sense that the resultant system collapses to its classical counterpart (Dunn 1988).\footnote{We also note a quite opposing viewpoint, discussed in Dunn (1980), which hearsay attributes to a lecture Saul Kripke gave at the University of Pittsburgh in 1974 (see Stairs (forthcoming)).  Kripke apparently argued that given a logicist or set-theoretic understanding of numbers, it can be shown using Putnam's views 
that $2 + 2 > 4$ since the cartesian product of a 2-membered set with itself has more than $4$ ordered pairs.}

In Dedekind's treatise on abstract structure of numbers, he asked \lq\lq what are numbers and what should they be?", then answered:\lq\lq numbers are free creations of the human mind.  They serve as a means of apprehending more easily and more sharply the difference of things."  Wave functions are creations of the human mind.  However, as a means of apprehending things, they are neither easier nor more sharply distinguished than numbers.  All's fair in love and quantum theory.

\vspace{.2in}
\textbf{References}
\vspace{.1in}

Bennett, C.H. and Brassard,G. (1984) "Quantum Cryptography: Public key distribution and coin tossing", in Proceedings of the IEEE International Conference on Computers, Systems, and Signal Processing, Bangalore, p. 175 (1984).

Birkhoff, G. (1940), \textit{Lattice Theory}, revised editions 1948 and 1967,
Colloquium Publications, vol. 25, American Mathemaical Society, Providence, RI.

Birkhoff, G. and von Neumann, J. (1936), \textquotedblleft The Logic of
Quantum Mechanics,\textquotedblright\ \textit{Annals of Mathematics}, 37, pp. 823--843.

Boole, G. (1847), \textit{The Mathematical Analysis of Logic, Being an Essay
Towards a Calculus of Deductive Reasoning}, Cambridge: Macmillan, Barclay, \&
Macmillan; reprinted Oxford: Basil Blackwell 1951.

Boole, G. (1854), \textit{An Investigation of The Laws of Thought on Which are
Founded the Mathematical Theories of Logic and Probabilities}, London:
Macmillan; reprint by Dover 1958.

Dalla Chiara, M L., Giuntini, R., and Greechie, R. J. (2004),
\textit{Reasoning in Quantum Theory: Sharp and Unsharp Quantum Logics}, Trends
in Logic - \textit{Studia Logica} Library, Kluwer Academic Publishers,

Dalla Chiara, M L. and Giuntini, R.\ \textquotedblleft Quantum Logics," in
\textit{Handbook of Philosophical Logic}, 2nd Edition, eds. D. M. Gabbay and
F. Guenthner, Kluwer Academic Publishers, Dordrecth, The Netherlands, pp. 129-228.

Dipert, R. (1978), \textit{Development and Crisis in Late Boolean Logic: The
Deductive Logics of Jevons, Peirce, and Schr\"{o}der}, Ph. D. Dissertation,
Indiana University.

Dunn, J. M., \textit{Quantum Mathematics}, Proceedings of the Biennial Meeting of the Philosophy of Science Association,
Vol. 1980, Volume Two: Symposia and Invited Papers (1980), pp. 512-531.

Dunn, J.M., \textit{The Impossibility of Certain Higher-Order Non-Classical Logics
with Extensionality}, David F. Austin, editor, Philosophical Analysis:
A Defense By Example (Dordrecht, Kluwer Academic
Publishers, 1988), xv + 363 pp.

Dunn, J. M., Moss, L., Hagge, T., and Wang, Z. (2005), \textquotedblleft
Quantum Logic as Motivated by Quantum Computing,\textquotedblright%
\ \textit{The Journal of Symbolic Logic} 70, pp. 353-359.

Egly, U. and Tompits, H. (1999), \textquotedblleft Gentzen-Like Methods in
Quantum Logic,\textquotedblright%
\ www.kr.tuwien.ac.at/staff/tompits/papers/tableaux-99.pdf, \ presented at
TAB-\newline LEAUX'99: International Conference on Analytic Tableaux and
Related Methods, Saratoga Springs, NY, June 1999.

Freedman, M.H., Kitaev, A., Larsen,M., and Wang,Z. (2003),
\emph{Topological quantum computation.} Bull. Amer. Math. Soc. (N.S.) 40 (2003), no. 1, 31--38.

Hagge, T. J. (2007),  $QL(\Bbb C^n)$ determines $n$. J. Symbolic Logic 72 (2007), no. 4, 1194--1196.

Hastings, M.B. (2009), Superadditivity of communication capacity using entangled inputs, Nature Physics 5, 255 - 257 (2009)

Huhn, A.P. (1972), Schwach distributive Verbande I, Acta Sci. Math. 33 (1972), 297-305.

Kleene, S. C. (1952), \textit{Introduction to Metamathematics,} D. van
Nostrand Company, Inc., Princeton, NJ.

Greechie, R. J. (1983), A Non-standard Quantum Logic with a strong Set of
States, in E. G. Beltrametti and B. C. van Fraassen (eds.), \textit{Current
Issues in Quantum Logic}, vol. 9 of Ettore Majorana International Science
Series, Plenum, NY, pp. 375-380.

Kneale, W. (1956), \textquotedblleft Boole and the Algebra of
Logic,\textquotedblright\ \textit{Notes and Records of the Royal Society of
London}, 12, pp. 53-63

Mackey, G. (1957), \textit{Mathematical Foundations of Quantum Theory},
Benjamin, New York.

Nishimura, H. (2009), \textquotedblleft Gentzen Methods in Quantum Logic,
\textquotedblright\ in \textit{Handbook of Quantum Logic and Quantum
Structures},. eds.K. Engesser, D. M. Gabbay, and D. Lehmann,
Elsevier/North-Holland, Amsterdam and New York.

Putnam, H. (1968), 'Is logic empirical?' in eds. R. Cohen and M. Wartofsky, {\it{Boston Studies in the Philosophy of Science}}, vol. 5 (Dordrecht: Reidel), pp. 216–--241. Reprinted as "The logic of quantum mechanics" in H. Putnam, {\it{Putnam, H. (1968), 'Is logic empirical?'}} in R. Cohen and M. Wartofsky
(eds.), Boston Studies in the Philosophy of Science, vol. 5, Reidel, Dordrecht, pp. 216–--241. Reprinted as "The Logic of Quantum Mechanics" in H. Putnam, {\it{Mathematics, Matter, and Method. Philosophical Papers}}, vol. 1. Cambridge University Press, Cambridge UK, 1975, pp. 174–--197.

R\'{e}dei, M. (2005), ed, \textit{John von Neumann : Selected Letters},
History of Mathematics, vol. 27, American Mathematical Society, Providence, RI.

R\'{e}dei, M. (2007), \textquotedblleft The Birth of Quantum
Logic\textquotedblright, \textit{History and Philosophy of Logic}, 28, 107--122

Schroeder, M. (1997), \textquotedblleft A Brief History of the Notation of
Boole's Algebra,\textquotedblright. \textit{Nordic Journal of Philosophical
Logic}, 2, 1997, pp. 41--62.

Shor, P. W. (1994), Algorithms for quantum computation: discrete logarithms and factoring. 35th Annual Symposium on Foundations of Computer Science (Santa Fe, NM, 1994), 124--134, IEEE Comput. Soc. Press, Los Alamitos, CA, 1994.

Stairs, A. (forthcoming), "Could Logic be Empirical? The Putnam Kripke Debate,". in Logic $\&$ Algebraic Structures in Quantum Computing $\&$ Information, eds. J. C. Reimann, V. Harizanov, and A. Eskandarian, Lecture Notes in Logic, Cambridge University Press, Cambridge UK.

M. H. Stone (1936), \textquotedblleft The Theory of the Representations of
Boolean Algebras,\textquotedblright\ \textit{Transactions of the American
Mathematical Society}, 40, 37-111.

von Neumann, J. (1932) \textit{Grundlagen der Quantenmechanik}, Springer
Verlag, Berlin, Heidelberg, New York. English translation,
\textit{Mathematical Foundations of Quantum Mechanics}, Princeton University
Press, Princeton, NJ, 1955.

Wang, Z. (2010), Topological quantum computation, CBMS monograph, vol 112,  Amer. Math. Soc., 2010.

\end{document}